\documentclass[reqno]{amsart}

\usepackage{graphicx}
\usepackage{amssymb}
\usepackage{yfonts}
\usepackage[all]{xy}
\usepackage[inline]{enumitem}

 \newcommand{\ev}{\mathop{\rm ev}\nolimits}

\begin{document}

\ \\[-30mm]
\title{Universal covering spaces, a footnote} 
\author[Petar Pave\v{s}i\'{c}]{Petar Pave\v si\'{c}}
\address{$^{*}$ Faculty of Mathematics and Physics, University of Ljubljana,
Jadranska 21,  1000 Ljubljana, Slovenija}
\email{petar.pavesic@fmf.uni-lj.si}


\begin{abstract}
We give a short proof that (for nice $X$) the based fundamental groupoid of $X$ with topology 
induced by the compact open topology on the space of paths, is indeed the universal covering space of $X$.
\ \\[3mm]
\end{abstract}
\maketitle
\ \\[-10mm]
Standard  textbook constructions of the universal covering space of $X$ are usually given as follows. 
Pick a point $x_0\in X$, consider the set $P(X,x_0)$ of all paths in $X$ starting at $x_0$ and identify
elements of $P(X,x_0)$ that are path-homotopic. (Recall that arbitrary paths $\alpha,\alpha'$ in $X$ 
are \emph{path-homotopic}, denoted $\alpha\sim\alpha'$, if $\alpha\simeq\alpha' (\mathrm{rel}\ \partial I)$.) 
Endow the set of equivalence classes $P(X,x_0)/\sim$ with the topology generated by subsets of the form 
$\langle \alpha,U\rangle$,
where $\alpha\in P(X,x_0)$, $U$ is an open neighbourhood of $\alpha(1)$ and 
$\langle \alpha,U\rangle$ consists of equivalence classes of paths obtained by concatenating $\alpha$ to a 
path in $P(U,\alpha(1))$.  It then 
requires certain amount of work (see for example \cite[pp. 82-83]{Spanier1966}) 
to show that the resulting space is the universal covering space of $X$
(provided $X$ satisfies suitable connectivity assumptions).

There is perhaps a more natural way to topologize $P(X,x_0)/\sim$, namely by taking the quotient topology 
of the standard compact-open topology on $P(X,x_0)$. It raises an obvious question: are the two 
topologies equal? A short answer is: not in general, but yes if $X$ admits a universal covering space.

So why are various textbooks not using quotient of the compact-open topology, which seems to be more standard
and easier to define? Indeed, Spanier \cite[p.82]{Spanier1966} states that 
\emph{`We could start with the compact-open topology on $P(X,x_0)$ and use the quotient topology on the 
set of equivalence classes, 
but it seems no simpler than merely topologizing the set of equivalence classes directly, as is done below.'}
One can have long discussions about more or less efficient approaches, but the fact is that most later authors 
followed Spanier's approach without further explanation which resulted in an unintended side-effect. 
I was asked on a number of occasions,
even by topologists, why is the compact-open topology inadequate for the construction of universal 
covering spaces. Similar questions pop-out on various forums so it appears that there is 
a widespread conviction that Spanier's approach is not just a matter of convenience but that compact-open 
topology is actually `wrong'.

Although one can find in the literature statements that the two topologies coincide (see notes at the end), they are
always buried among more general results.
That's why I decided to write a short direct proof that for nice spaces universal coverings can be constructed 
by taking 
quotient topology of the compact-open topology on the space of based paths. The proof is reasonably 
self-contained: in the next Lemma we prove some properties of the path-space fibration using only 
standard facts from general topology (as in Dugundji \cite{Dugundji1966}); the proof of the main  
result then easily follows, using only the unique lifting property of  covering spaces.

For an arbitrary space $X$ and $x_0\in X$ let $P(X,x_0)$ denote the space of based continuous paths 
$\alpha\colon (I,0)\to (X,x_0)$, endowed with the compact-open topology and let 
$\ev\colon P(X,x_0)\to X$ denote the evaluation map given as $\ev(\alpha):=\alpha(1)$. 
It is well-known that $\ev$ is continuous (it is even a Hurewicz fibration, but that fact will not be needed). 
Before stating the relevant properties of evaluation maps, let us recall the definition of 
semi-locally simply connected spaces. 
We will say that $U\subseteq X$ is a \emph{slsc}-subset of $X$ if every loop in $U$ can be contracted 
to a point in $X$ (or equivalently, if every two paths in $U$ with common endpoints
are path-homotopic in $X$, or that inclusion induces a trivial homomorphism
$\pi_1(U,x)\to \pi_1(X,x)$ for every $x\in U$). A space $X$ is \emph{semi-locally simply-connected} 
if it can be covered by open slsc-subsets.

\textbf{LEMMA}\\[-7mm]
\begin{enumerate}
\item If $X$ is path-connected, then $\ev$ is surjective.
\item If $X$ is locally path-connected, then $\ev$ is an open map.
\item If $X$ is locally path-connected and semi-locally simply-connected, 
then the path-components of $\ev^{-1}(U)$ are open in $P(X,x_0)$ for every open subset $U\subseteq X$.
\end{enumerate}
\begin{proof}
(1) is trivial. 
In the proofs of (2) and (3) we will follow Dugundji \cite[Ch. XII]{Dugundji1966} and denote by 
$(A,V):=\{f\colon I\to X\mid f(A)\subset V\}$ ($A\subseteq I$ compact, $V\subseteq X$ open) the sub-basic
open sets that generate compact-open topology on $P(X,x_0)$. 

(2) Consider a point in $\ev\big(\bigcap_1^n (A_i,V_i)\big)$ of the form $\alpha(1)$ for some 
$\alpha\in \bigcap_1^n (A_i,V_i)$ and let $V$ be the path-component 
of $\bigcap\{V_i\mid  1\in A_i\}$ containing $\alpha(1)$.  Since $X$ is locally path-connected,
$V$ is open in $X$. For any given point $v\in V$ one can extend the final segment of $\alpha$ 
and obtain a path in $\bigcap_1^n (A_i,V_i)$ ending at $v$ (cf. \cite[II,4.4]{Spanier1966}), 
therefore $\alpha(1)$ is an interior point of $\ev\big(\bigcap_1^n (A_i,V_i)\big)$.

(3) Given $\alpha\in \ev^{-1}(U)$ one can find points $0=t_0<t_1<\cdots<t_n=1$ and open slsc-sets 
$U_1,\ldots,U_n$, such that $\alpha([t_{i-1},t_i])\subset U_i$. Moreover, one can find 
path-connected open sets $V_1,\ldots,V_n$, such that $\alpha(t_i)\in V_i\subseteq U_i\cap U_{i+1}$ 
for $i=1,\ldots,n-1$ and $\alpha(1)\in V_n\subseteq U_n\cap U$. Then 
$W:=\bigcap_{i=1}^n ([t_{i-1},t_i],U_i) \cap \bigcap_{i=1}^n (\{t_i\},V_i)$ is an open neighbourhood 
of $\alpha$ in $P(X,x_0)$. To show that $W$ is contained in a path-component of $\ev^{-1}(U)$ we 
must verify that every element of $W$ is connected to $\alpha$ by a path in $\ev^{-1}(U)$. 
Note that a path from $\alpha$ to $\beta$ in $\ev^{-1}(U)$ corresponds to a path-homotopy 
$\alpha\cdot\lambda\sim \beta$ for some $\lambda\colon(I,0,1)\to (U,\alpha(1),\beta(1))$. 
Given $\beta\in W$ choose paths $\lambda_i\colon (I,0,1)\to (V_i,\alpha(t_i),\beta(t_i))$
for $i=1,\ldots,n$. The paths $\overline\lambda_{i-1}\cdot\alpha|_{[t_{i-1},t_i]}\cdot\lambda_i$ and 
$\beta|_{[t_{i-1},t_i]}$ are contained in $U_i$ and are thus path-homotopic. By concatenation we obtain 
a path-homotopy $\alpha\cdot\lambda_n\sim \beta$, so $\alpha$ and $\beta$ 
are connected by a path in $\ev^{-1}(U)$. Therefore $\alpha$ is an interior point of a 
path-component of $\ev^{-1}(U)$ as claimed.
\end{proof}

We are now ready to construct the universal covering for $X$. 
Let \emph{(topological) fundamental groupoid of $X$ based at $x_0$} be the quotient topological space
$\Pi(X,x_0):=P(X,x_0)/\!\!\sim$ 
with the projection to equivalence classes $q\colon P(X,x_0)\to \Pi(X,x_0)$ as a quotient map.
Since $\ev$ preserves $\sim$, we obtain by \cite[VI,4.3]{Dugundji1966} the induced continuous map 
$p\colon\Pi(X,x_0)\to X$, $p([\alpha]):=\alpha(1)$ for which $p\circ q=\ev$.

\textbf{THEOREM} \emph{If $X$ is path-connected, locally path-connected and semi-locally simply connected, 
then $p\colon\Pi(X,x_0)\to X$ is the universal covering of $X$.}\\[-7mm]
\begin{proof}
To begin with, $p$ is in open map by \cite[VI, 3.2(2)]{Dugundji1966} and (2) of the above Lemma.
We claim that every path-connected open slsc-subset $U$ of $X$ is evenly covered by $p$
in the sense of \cite[II.1]{Spanier1966}. By Lemma (2) all path-components of $\ev^{-1}(U)$ are open, therefore 
$\ev^{-1}(U)$ is homeomorphic to the topological
disjoint union of its path-components. As the relation $\sim$ identifies points within the same 
path-component, to 
compute $p^{-1}(U)=\ev^{-1}(U)/\!\sim$ it is sufficient to examine the quotients of individual path-components. 
For each path-component $C$ of $\ev^{-1}(U)$ consider the restriction $p_C\colon (C/\!\!\sim)\to U$.  
Since $U$ is path-connected, $p_C$ is clearly surjective. To check injectivity, take $\alpha,\alpha'\in C$, 
which means that there exists a path $\lambda\colon(I,0,1)\to (U,\alpha(1),\alpha'(1))$, such that 
$\alpha'\sim\alpha\cdot\lambda$. If $p_C([\alpha])=p_C([\alpha'])$ (equivalently, if $\alpha(1)=\alpha'(1)$),
then $\lambda$ is a loop in the slsc-set $U$, therefore $[\alpha]=[\alpha']$. Finally, $p_C$ is  
continuous and open, because it is a restriction of a continuous open map $p$ to the open
subspace $C/\!\!\sim$ of $\Pi(X,x_0)$. Therefore $p^{-1}(U)$ is the disjoint union of open subsets, each 
of which is mapped homeomorphically to $U$ by $p$.

Every path $\alpha\in P(X,x_0)$ can be lifted to a path in $\widetilde\alpha\colon I\to \Pi(X,x_0)$, given
explicitly by the formula $\widetilde\alpha(t):=[s\mapsto \alpha(st)]$. Clearly, $\widetilde\alpha(0)=c_{x_0}$,
$\widetilde\alpha(1)=\alpha$ and the lifting is unique, because $p$ is a covering projection.
Recall \cite[II,3]{Spanier1966} that the fundamental group of $\Pi(X,x_0)$ can be identified 
with the subgroup of $\pi_1(X.x_0)$ represented by loops whose lifting along $p$ is also a loop.
But given a based loop $\alpha$ in $(X,x_0)$ its unique lifting $\widetilde\alpha$ is a loop if, and only 
if $\alpha\sim $, therefore $\Pi(X,x_0)$ is simply-connected, and thus universal by 
\cite[II,5.7]{Spanier1966}.
\end{proof}

\newpage

{\bf Sub-footnotes:}

(1) Spanier's book is a standard reference for the construction of universal coverings, but the 
original idea is much older. I was able to date it at least to H. Weyl's monograph on 
Riemann surfaces from 1913. Here is an excerpt from a 1955 translation
with an extraordinary modern description of the universal covering surface.  

\hspace{5mm}\begin{minipage}{140mm}
\emph{The universal covering surface may be defined as follows. If $\textfrak{p}_0$ is a fixed point 
of \textfrak{F} then every curve $\gamma$ starting at $\textfrak{p}_0$ defines a ``point of 
$\tilde{\textfrak{F}}$'' which we say lies over the endpoint of $\gamma$. Two such curves $\gamma,\gamma'$ define 
the same point of $\tilde{\textfrak{F}}$ if and only if on every perfect covering surface over \textfrak{F} every 
pair of curves which start at the same
point have traces $\gamma,\gamma'$ end at the same point. Let $\gamma_0$ be a curve on \textfrak{F} from 
$\textfrak{p}_0$ to $\textfrak{p}$ which defines the point $\widetilde{\textfrak{p}}$  on 
$\tilde{\textfrak{F}}$, and let \textfrak{U} be a neighbourhood of $\textfrak{p}$ on 
\textfrak{F}. If I attach to $\gamma_0$ all possible curves $\gamma$ in \textfrak{U} 
which start at $\textfrak{p}$, then I say that the points of $\tilde{\textfrak{F}}$ defined by all these 
curves $\gamma_0+\gamma$ form a ``neighborhood'' $\tilde{\textfrak{U}}$
of $\tilde{\textfrak{p}}$. Since \textfrak{U} is simply connected there is just one point of 
$\tilde{\textfrak{U}}$ over each point of \textfrak{U}; hence our concept of 
``neighborhood'' satisfies the conditions stated in \S 4. } 
(see \cite[pp. 58-59]{Weyl1913})
\end{minipage}

(2) Various topologies on the based fundamental groupoid were thoroughly studied in the recent years. 
In that context the topology used by Spanier was called \emph{whisker topology}, while the quotient of 
the compact-open topology was called \emph{CO-topology} (see Virk and Zastrow \cite{VirkZastrow2014}). In particular, 
Fisher and Zastrow \cite[Lemma 2.1]{FischerZastrow2007} proved that whisker topology is always stronger 
than CO-topology and that the two agree for spaces that are locally path-connected and semi-locally simply 
connected.
However, there are locally path-connected spaces (e.g., the `Hawaiian earring') for which whisker topology is 
strictly stronger than CO-topology.

(3) It can be occasionally useful to know that the fact that $\ev \colon P(X,x_0)\to X$ is an open map, 
that played an important role in our argument, can be extended to arbitrary Hurewicz fibrations. Indeed,
McAuley \cite[Cor.\,1.1]{McAuley1985} proved that $X$ is locally path-connected if, and only if, every 
Hurewicz fibration with base $X$ is an open map.

\section*{Acknowledgment}

I am indebted with Stefan Friedl who kindly pointed out that there was a gap in the proof in a previous 
version of this note and suggested how it could be filled. Further discussion lead to a considerably shorter 
and simpler proof of the main result (see \cite[Part IV, Sec.\,29]{Friedl} for more information).

\end{document}